\documentclass[11pt]{amsart}
\numberwithin{equation}{section}
\newtheorem{lemma}{Lemma}

\numberwithin{theorem}{section} \numberwithin{lemma}{section}
\numberwithin{proposition}{section} \numberwithin{equation}{section}

\def\al{\aligned}
\def\eal{\endaligned}

\begin{document}

\tracingpages 1
\title[Ricci flow with surgery]{\bf Strong non-collapsing and uniform Sobolev inequalities for Ricci flow with surgeries}
\author{ Qi S. Zhang}

\address{Department of Mathematics,  University of California,
Riverside, CA 92521, USA }
\date{December 2007}

\maketitle



\begin{abstract}

We prove a uniform Sobolev inequality for Ricci  flow, which is
independent of the number of surgeries.   As an application, under
less assumptions, a non-collapsing result stronger than Perelman's
$\kappa$ non-collapsing with surgery is derived. The proof is
shorter and seems more accessible. The result also improves some
earlier ones where the Sobolev inequality depended on the number of
surgeries.

\end{abstract}


\section{Introduction}


A crucial step in Perelman's work on Poincar\'e and Geometrization
conjectures is the $\kappa$ non-collapsing result for Ricci flow
with or without surgeries. The proof of this result in the surgery
case requires truely complicated calculation using such new concepts
as reduced distance, admissible curve, barely admissible curve,
gradient estimate of scalar curvature etc. This is elucidated in
great length by Cao and Zhu [CZ], Kleiner and Lott [KL] and Morgan
and Tian [MT].

In this paper we prove a uniform Sobolev inequality for Ricci  flow,
which is independent of the number of surgeries.  It is well known
that uniform Sobolev inequalities are essential in that they encode
rich analytical and geometrical information on the manifold. These
include, non-collapsing, isoperimetric inequalities etc. As a
consequence, a strong non-collapsing result is obtained. It includes
Perelman's $\kappa$ non-collapsing with
 surgery as a special case.  The result also requires less assumptions.
 For instance we do not need the canonical neighborhood assumption for the
 whole manifold(see Remark 1.2 below).
  In the proof, we use only Perelman's
$W$ entropy and some analysis of the minimizer equation of the $W$
entropy on horn like manifolds. Hence it is shorter and seems more accessible.

Let ${\bf M}$ be a compact Riemannian manifold of dimension $n \ge
3$ and $g$ be the metric. Then a Sobolev
inequality of the following form holds: there exist positive
constants $A, B$
 such that, for all $v
\in W^{1, 2}({\bf M}, g)$,
\begin{equation}
\label{sob}
 \bigg( \int v^{2n/(n-2)} d\mu(g) \bigg)^{(n-2)/n} \le
A \int |\nabla v |^2 d\mu(g) + B \int v^2 d\mu(g).
\end{equation}

This inequality was proven by Aubin \cite{Au:1}  for
$A=K^2(n)+\epsilon$ with $\epsilon>0$ and $B$ depending on bounds on
the injectivity radius, sectional curvatures. Here $K(n)$ is the
best constant in the Sobolev imbedding for ${\bf R}^n$. Hebey
\cite{H:1} showed that $B$ can be chosen to depend only on
$\epsilon$, the injectivity radius and the lower bound of the Ricci
curvature. Hebey and Vaugon [HV] proved that one can even take
$\epsilon=0$. However the constant $B$ will also depend on the
derivatives of the curvature tensor. Hence, the controlling
geometric quantities for $B$ as stated above are not invariant under
the Ricci flow in general.  Theorem 1.1 below states that a uniform
Sobolev inequality of the above type holds uniformly under Ricci
flow in finite time, even in the presence of indefinite number of
surgeries.

In order to state the theorem, we first introduce some notations. They are mainly taken
from [P1,2], [CZ], [KL] and [MT].

We use $({\bf M}, g(t))$ to denote Hamilton's Ricci flow, $\frac{d
g}{d t}= - 2 Ric$. If a surgery occurs at time $t$, then $({\bf M},
g(t^-))$ denotes the pre surgery manifold (the one right before the
surgery); and $({\bf M}, g(t^+))$ denotes the post surgery manifold
 (the one right after the surgery).
  The ball of radius $r$ with respect to the metric $g(t)$, centered at $x$,
  is denoted by $B(x, t, r)$.  The scalar curvature is denoted by
  $R=R(x, t)$ and $R^-_0 = \sup R^-(x, 0)$. $Rm$ denotes the full curvature tensor.
  $d\mu(g(t))$ denotes the volume element. $vol({\bf M(g(t))}$ is the total volume
  of ${\bf M}$ under $g(t)$.

 In this paper we use the following definition of $\kappa$ non-collapsing by Perelman [P2],  as elucidated in Definition 77.9 of [KL].

{\it {\bf Definition 1.1.} $\kappa$ non-collapsing.

Let $({\bf M}, g(t))$ be a Ricci flow with surgery defined on $[a,
b]$. Suppose that $x_0 \in {\bf M}$, $t_0 \in [a, b]$ and $r>0$ are
such that $t_0-r^2 \ge a$, $B(x_0, t_0, r) \subset {\bf M}$ is a
proper ball and the parabolic ball $P(x_0, t_0, r, -r^2)$ is
unscathed. Then ${\bf M}$ is $\kappa$-collapsed at $(x_0, t_0)$ at
scale $r$ if $|Rm| \le r^{-2}$ on $P(x_0, t_0, r, -r^2)$ and
$vol(B(x_0, t_0, r)) < \kappa r^3$; otherwise it is
$\kappa$-noncollapsed.}

\medskip
Here we introduce
\medskip

{\it {\bf Definition 1.2.} Strong $\kappa$ non-collapsing.

Let ${\bf M}$ be a Ricci flow with surgery defined on $[a, b]$.
Suppose that $x_0 \in {\bf M}$, $t_0 \in [a, b]$ and $r>0$ are such
that $B(x_0, t_0, r) \subset {\bf M}$ is a proper ball . Then ${\bf
M}$ is strong $\kappa$-noncollapsed at $(x_0, t_0)$ at scale $r$ if
$R \le r^{-2}$ on $B(x_0, t_0, r)$ and $vol(B(x_0, t_0, r)) \ge
\kappa r^3$. }
\medskip

This strong $\kappa$ non-collapsing improves the $\kappa$ non-collapsing on two aspects. One is that only information on the metric balls on one time level is needed.
Thus it bypasses the complicated issue that a parabolic ball may be cut by a surgery.
The other is that it only requires scalar curvature upper bound.
\medskip

{\it {\bf Definition 1.3.}   Normalized manifold.

 A compact  Riemannian manifold is normalized if $|Rm| \le 1$ everywhere and the volume of every unit ball is at least half of the volume of the Euclidean unit ball.}
\medskip

{\it {\bf Definition 1.4.}  $\epsilon$ neck,  $\epsilon$ horn,
double $\epsilon$ horn, and $\epsilon$ tube.

 An $\epsilon$ neck (of radius $r$) is an open set with a metric which is, after scaling the metric with factor $r^{-2}$, $\epsilon$ close, in the
$C^{\epsilon^{-1}}$ topology,  to the standard neck $S^2 \times
(-\epsilon^{-1}, \epsilon^{-1})$. Here and later $C^{\epsilon^{-1}}$
means $C^{[\epsilon^{-1}]+1}$.

Let $I$ be an open interval in $R^1$. An $\epsilon$ horn (of radius $r$)  is  $S^2 \times I$ with a metric with the following properties:  each point is
contained in some $\epsilon$ neck;   one end is contained in an $\epsilon$ neck of radius
$r$; the scalar curvature tends to infinity at the other end.}

An $\epsilon$ tube  is  $S^2 \times I$ with a metric such that  each point is
contained in some $\epsilon$ neck and the scalar curvature stays bounded on both
ends.

A double $\epsilon$ horn  is  $S^2 \times I$ with a metric such that
each point is contained in some $\epsilon$ neck and the scalar
curvature tends to infinity at both ends.
\medskip

{\it {\bf Definition 1.5.}  A standard capped infinite cylinder is
${\bf R^3}$ equipped with a rotationally symmetric metric with
nonnegative sectional curvature and positive scalar curvature such
that outside a compact set it is a semi-infinite standard round
cylinder $S^2 \times (-\infty, 0)$.}

\medskip

A few more basic facts concerning Ricci flow with surgery such as $(r, \delta)$ surgery, $\delta$ neck is given in the appendix.   For detailed information and related
terminology on
Ricci flow with surgery we refer the reader to [CZ], [KL] and [MT].

Here is the main result of paper.
\medskip

{\it {\bf Theorem 1.1.} Given real numbers $T_1<T_2$, let $({\bf
M}, g(t))$ be a 3 dimensional Ricci flow with normalized initial
condition defined on the time interval containing $[T_1, T_2]$.
Suppose the following conditions are met.


(a). There are finitely many $(r, \delta)$ surgeries in $[T_1,
T_2]$, occurring in $\epsilon$ horns of radii $r$. Here $r \le r_0$
and $\epsilon \le \epsilon_0$, with $r_0$ and $\epsilon_0$ being
fixed sufficiently small positive numbers less than $1$. The surgery radii are $h
\le \delta^2 r$ i.e. the surgeries occur in $\delta$ necks of radius
$h \le \delta^2 r$. Here $0<\delta \le \delta_0$ where
$\delta_0=\delta_0(r, \epsilon_0)>0$ is sufficiently small. Outside
of the $\epsilon$ horns, the Ricci flow is smooth.

(b).  
For a constant $c>0$ and any point $x$ in all the above $\epsilon$
horns, the following holds:
there is a region $U$, satisfying, $B(x,  c \epsilon^{-1} R^{-1/2}(x)) \subset U \subset B(x, 2 c \epsilon^{-1} R^{-1/2}(x))$, such that,  after scaling by a factor $R(x)$, it is 
$\epsilon$ close in the $C^{\epsilon^{-1}}$ topology to 
$S^2 \times (-\epsilon^{-1}, \epsilon^{-1})$.

Also for any $x$ in the  modified part of the $\epsilon$
horn immediately after a surgery, the following holds: the ball
$B(x,  \epsilon^{-1} R^{-1/2}(x))$, is, after scaling by a factor $R(x)$,
$\epsilon$ close in the $C^{\epsilon^{-1}}$ topology to the
corresponding ball of the standard capped infinite cylinder.

(c). For  $A_1>0$, the Sobolev imbedding
\[
\bigg( \int v^{2n/(n-2)} d\mu(g(T_1)) \bigg)^{(n-2)/n} \le A_1 \int (
4 |\nabla v |^2 +  R v^2 ) d\mu(g(T_1)) + A_1  \int v^2
d\mu(g(T_1)). \
\]holds for all $v \in W^{1, 2}({\bf M}, g(T_1))$.

Then for all $t \in (T_1, T_2]$, the Sobolev imbedding below holds
for all $v \in W^{1, 2}({\bf M}, g(t))$.
\[
\bigg( \int v^{2n/(n-2)} d\mu(g(t)) \bigg)^{(n-2)/n} \le A_2 \int (
4 |\nabla v |^2 + R v^2 ) d\mu(g(t)) + A_2  \int v^2
d\mu(g(t))
\]
Here
\[
A_2   = C(A_1, \sup R^-(x, 0), T_2, T_1, \sup_{t \in [T_1, T_2]}
Vol({\bf M}(g(t))) \ )
\] is independent of the number of surgeries or $r$.

Moreover, the Ricci flow is strong $\kappa$ noncollapsed in the whole
interval $[T_1, T_2]$  under  scale $1$ where $\kappa$ depends only on $A_2$.
 }

\medskip

{\it Remark 1.1}. By the work  Hebey \cite{H:1}, at any given time,
a Sobolev imbedding always holds with constants depending on lower
bound of Ricci curvature and injectivity radius. So one can replace
assumption (c) by the assumption that $({\bf M}, g(T_1))$ is
$\kappa$ non-collapsed and that the canonical neighborhood
assumption (with a fixed radius $r_0>0$ and $\epsilon_0>0$) at $T_1$
holds. It is easy to see that these together imply the Sobolev
imbedding at time $T_1$.

We assume as usual that, at a surgery, we throw away all compact
components with positive sectional curvature , and also capped
horns, double horns and all compact components lying in the region
where $R>\delta r$. In the extra assumption that the Ricci flow is
smooth outside of the $\epsilon$ horns, we have excluded these
deleted items.

{\it Remark 1.2}. With the exception of using the monotonicity of
Perelman's $W$ entropy, the proof of Theorem 1.1 uses only long
established results.  Under $(r, \delta)$ surgery,  assumption (b) is 
clearly implied by, but
much weaker than the {\it canonical neighborhood assumption} on the
whole manifold ${\bf M}$, which was used in all the papers so far.

{\it Remark 1.3}  In [Z2], it was shown that under a Ricci flow with
finite number of surgeries in finite time, a uniform Sobolev
imbedding holds.  In [Y], a similar result depending on the number
of surgeries was stated without proof.

Let us finish the introduction by outlining the proof.  Recall Perelman's
$W$ entropy and its monotonicity.  They are in fact the monotonicity of the best
constants of the Log Sobolev inequality with certain parameters.  If a Ricci flow
is smooth over a finite time interval, then the best constants of the Log Sobolev
inequality with a changing parameter does not decrease.
If a Ricci flow undergoes a $(r, \delta)$ surgery with $\delta$ sufficiently small,
then the best constant only decreases by at most a constant times the change
in volume.  This is achieved by a weighted estimate of Agmon type for the
minimizing equation of the $W$ entropy. The method is motivated by those
at the end of [P2] and [KL] where the change of eigenvalues of the linear operator
$4 \Delta -R$ was studied.  Therefore  in  finite time, the best constant
of the Log Sobolev inequality with certain parameters is uniformly bounded from below by a negative constant, regardless of the number of surgeries.
This uniform Log Sobolev inequality is then converted by known method to the desired uniform
Sobolev inequality which in turn yields strong noncollapsing.

\section{Proof of Theorem 1.1.}

We will need three lemmas before carrying out the proof of the theorem.
Much of the analysis is focused on the $\epsilon$ horn where a surgery takes place.
So we will fix some notations concerning the $\epsilon$ horn and the surgery cap.

Recall that a $(r, \delta)$ surgery occurs deep inside a $\epsilon$
horn of radius $r$. The horn is cut open at the place where the
radius is $h \le \delta^2 r$. Then a cap is attached and a smooth
metric is constructed by interpolating between the metric on the
horn and the metric on the cap. The resulting manifold right after
the surgery is denoted by ${\bf M^+}$ and the $\epsilon$ horn thus
modified by the surgery is called a capped $\epsilon$ horn with
radius $r$.

Let $D$ be a capped $\epsilon$ horn.  By assumption, a region ${\bf
N}$ around the boundary $\partial D$ equipped with the scaled metric
$c r^{-2} g$ is $\epsilon$ close, in the $C^{\epsilon^{-1}}$
topology, to the standard round neck $S^2 \times (-\epsilon,
\epsilon)$. Here $c$ is a generic positive constant such that $c
r^{-2}$ equals the scalar curvature at a point on $\partial D$. For
this reason we will often take $c=1$.

Let $\Pi$ be the diffeomorphism from the
standard round neck to ${\bf N}$ in the definition of $\epsilon$
closeness. Denote by $z$ for a number in $(-\epsilon^{-1},
\epsilon^{-1})$.  For $\theta \in S^2$, $(\theta, z)$ is a
parametrization of ${\bf N}$ via the diffeomorphism $\Pi$. We can
identify the metric on ${\bf N}$ with its pull back on the round
neck by $\Pi$ in this manner.  We normalize the parameters so that the capped
$\epsilon$ horn lies in the region where $z \ge 0$

Next we define
\[
Y(D) = \inf \{  \frac{ 4 \int |\nabla v |^2 + R v^2}{\left(\int v^{2n/(n-2)} \right)^{(n-2)/n} }\
|\  v \in C^{\infty}_0(D \cup {\bf N}), v>0 \}.
\leqno(2.0)
\]

This quantity is clearly a conformal invariant. Since $R$ is positive in $D \cup {\bf N}$,
Y(D) is bounded from above and below by constant multiples of the Yamabe constant
\[
Y_0(D) =  \inf \{  \frac{ 4 \frac{n-1}{n-2} \int |\nabla v |^2 + R
v^2}{\left(\int v^{2n/(n-2)} \right)^{(n-2)/n} }\ |\  v \in
C^{\infty}_0(D \cup {\bf N}), v>0 \}.
\]

Let $g=g(x)$ be the metric on $D \cup {\bf N}$ then $Y(D)$ and
$Y_0(D)$ stay the same under the metric $g_1(x)= R(x) g(x)$.

Consider the manifold $(D \cup {\bf N}, g_1)$. By assumption and the
$(r, \delta)$ surgery procedure, there is a fixed $r_0>0$ such that
for any $x \in D \cup {\bf N}$, the ball $B(x, r_0)$ under $g_1$ is
$\epsilon$ close (in $C^{\epsilon^{-1}}$ topology) to a part of the
standard capped infinite cylinder. Since $\epsilon$ is sufficiently
small, we know that the injectivity radius of $(D \cup {\bf N},
g_1)$ is bounded from below by a positive constant; and its Ricci
curvature is bounded from below. Actually it is easy to see that
these hold for a much larger domain containing $(D \cup {\bf N},
g_1)$.  By Proposition 6 in [H1], we can find a positive constant
$C$ such that
\[
 \bigg( \int v^{2n/(n-2)} d\mu(g_1) \bigg)^{(n-2)/n} \le
C \int |\nabla v |^2 d\mu(g_1) + C \int v^2 d\mu(g)
\]for all $v \in
C^\infty_0(D \cup {\bf N})$. Since the scalar curvature of $(D \cup
{\bf N}, g_1)$ is bounded between two positive constants, we have,
for a constant still named $C$,
\[
 \bigg( \int v^{2n/(n-2)} d\mu(g_1) \bigg)^{(n-2)/n} \le
C \int (4 \frac{n-1}{n-2}|\nabla v |^2 + R v^2) d\mu(g_1)
\]for all $v \in
C^\infty_0(D \cup {\bf N})$.

Hence we see that $ Y_0(D) $ is bounded from below by a positive
constant when $\epsilon$ is sufficiently small. It is also bounded
from above by the Yamabe constant of $S^n$. Since $Y_0(D)$ and $Y(D)$ are
comparable,  we have shown
that
\[
0< Const_1 \le Y(D) \le Const_2 \leqno(2.1)
\]when $\epsilon$ is sufficiently small.

Next we present

\begin{lemma}
Let $({\bf M^+}, g)$ be a manifold right after a $(r, \delta)$
surgery.
Let $D \subset {\bf M^+} $ be a capped $\epsilon$ horn of radius $r$.
Here $\epsilon$ is a sufficiently small positive number.

Suppose $u$  with $\Vert u \Vert_{L^2(\bf M^+)}=1$ is a positive
solution to the equation
\[
\sigma^2 ( 4 \Delta u - R u ) + 2 u \ln u + \Lambda u + n( \ln
\sigma) u =0. \leqno(2.2)
\]
Here $\sigma>0$ and $\Lambda \le 0$.

Then there exists a positive constant $C$ depending only on $Y(D)$,
$n$ but not on the smallness of $\epsilon$ such that
\[
\sup_D u^2 \le C \max ( r^{-n}, \sigma^{-n} ).
\]
\proof
\end{lemma}

After taking the scaling
\[
g_1= \sigma^{-2} g, \ R_1 = \sigma^2 R, \ u_1 = \sigma^{n/2} u
\]we see that $u_1$ satisfies
\[
4 \Delta_1 u _1 - R_1 u_1 + 2 u_1 \ln u_1 + \Lambda u_1 =0.
\]Since the result in the lemma is independent of the above scaling,
we can just prove it for $\sigma = 1$.

So let $u$ be a positive solution to the equation
\[
 4 \Delta u - R u  + 2 u \ln u + \Lambda u  =0
\] in ${\bf M^+}$ such that its $L^2$ norm is 1.
Given any $p \ge 1$, it is easy to see that
\[
- 4 \Delta u^p + p R u^p \le 2 p u^p \ln u.
\leqno(2.3)
\]

We select a smooth cut off function $\phi$ which is one in $D$ and $0$ outside of
$D \cup {\bf N}$. Writing
$w = u^p$ and using $w \phi^2$ as a test function in (2.3), we deduce
\[
4 \int \nabla (w \phi^2) \nabla w + p \int R (w \phi)^2 \le 2 p \int (w \phi)^2 \ln u.
\]Since the scalar curvature $R$ is positive in the support of $\phi$ and $p \ge 1$, this shows
\[
4 \int \nabla (w \phi^2) \nabla w + \int R (w \phi)^2 \le  p \int (w \phi)^2 \ln u^2.
\]Using integration by parts, we have
\[
4 \int |\nabla (w \phi)|^2 + \int R (w \phi)^2 \le  4 \int |\nabla \phi|^2 w^2 +
p \int (w \phi)^2 \ln u^2.
\leqno(2.4)
\]

We need to dominate the  last term in (2.3) by the left hand side of (2.3).
For one positive number $a$ to be chosen later, it is clear that
\[
\ln u^2 \le u^{2a} + c(a).
\]Hence for any fixed $q>n/2$, the H\"older inequality implies
\[
\begin{aligned}
p \int (w \phi)^2 \ln u^2 &\le p \int (w \phi)^2 u^{2a} + p c(a) \int (w \phi)^2\\
&\le p \left( \int u^{2 a q} \right)^{1/q} \   \left( \int (w \phi)^{2  q/(q-1)} \right)^{(q-1)/q}
+ p c(a) \int (w \phi)^2.
\end{aligned}
\]We take $a = 1/q$ so that $2 a q =2$.
Since the $L^2$ norm of $u$ is $1$ by assumption, the above implies
\[
p \int (w \phi)^2 \ln u^2  \le  p  \left( \int (w \phi)^{2  q/(q-1)}
\right)^{(q-1)/q}
+ p c(a) \int (w \phi)^2.
\]By interpolation inequality (see p84 [HL] e.g.), it holds, for any $b>0$,
\[
\left( \int (w \phi)^{2  q/(q-1)} \right)^{(q-1)/q}
\le b
\left( \int (w \phi)^{2 n/(n-2)} \right)^{(n-2)/n}  + c(n, q) b^{-n/(2q-n)}
 \int (w \phi)^2.
\]Therefore
\[
p \int (w \phi)^2 \ln u^2  \le  p b \left( \int (w \phi)^{2 n/(n-2)}
\right)^{(n-2)/n}  + c(n, q) p b^{-n/(2q-n)} \int (w \phi)^2
+ p c(a) \int (w \phi)^2.
\leqno(2.5)
\]

By the definition of $Y(D)$ in (2.0),
we see that (2.4) gives

\[
Y(D) \left(\int w^{2n/(n-2)} \right)^{(n-2)/n} \le 4 \int |\nabla \phi|^2 w^2 +
p \int (w \phi)^2 \ln u^2.
\leqno(2.6)
\]Substituting (2.5) to the right hand side of (2.6), we arrive at
\[
\al
& Y(D) \left(\int w^{2n/(n-2)} \right)^{(n-2)/n} \le 4 \int |\nabla \phi|^2 w^2\\
 &\qquad +
 p b \left( \int (w \phi)^{2 n/(n-2)} \right)^{(n-2)/n}  + c(n, q) p b^{-n/(2q-n)} \int (w \phi)^2
+ p c(a) \int (w \phi)^2.
\eal
\]Take $b$ so that $p b = Y(D)/2$. It is clear that exist positive constant $c =c(Y(D), n, q)$
and $\alpha=\alpha(n, q)$ such that
\[
\left(\int w^{2n/(n-2)} \right)^{(n-2)/n} \le c (p+1)^{\alpha} \int (|\nabla \phi|^2+1) w^2.
\leqno(2.7)
\]

From here one can use Moser's iteration to prove the desired bound.
Let $z$ be the longitudinal parameter for $D$ described before the
lemma. For $z_2$ and $z_1$ such that $-1 \le z_2 < z_1 <0$, we
construct a smooth function of $z$, called $\xi$ such that
$\xi(z)=1$ when $z \ge z_1$; \ $\xi(z)=0$ when $z<z_2$ and $\xi(z)
\in (0, 1)$ for the rest of $z$. Set the test function $\phi=
\xi(z)=\xi(z(x))$. Then it is clear that
\[
|\nabla \phi| \le \frac{c}{r (z_1-z_2)}.
\leqno(2.8)
\]Write
\[
D_i=\{x \in {\bf M^+} \ | \ z(x) > z_i \}, \qquad i=1, 2.
\]By (2.7) and (2.8)
\[
\left(\int_{D_1} w^{2n/(n-2)} \right)^{(n-2)/n} \le c \max \{ \frac{1}{[(z_1-z_2) r]^2} , 1 \}
 (p+1)^{\alpha} \int_{D_2} w^2.
\leqno(2.9)
\]

Recall that $w = u^p$. We iterate (2.9) with $p=(n/(n-2))^i$, $i=0, 1, 2, ...$ in conjunction
with choosing
\[
z_1=- (1/2 + 1/2^{i+2}), \qquad z_2 = - (1/2 + 1/2^{i+1}).
\]Following Moser, we will get
\[
\sup_D u^2 \le C \max (r^{-n}, 1) \int u^2.
\]\qed

The next lemma is a nonlinear version of the result in [P2] and Lemma 92.10 in [KL].
This estimate has its origin in the weighted Agmon type estimate of eigenfunctions of the Laplacian.
\medskip

\begin{lemma}
Let $({\bf M}, g)$ be any compact manifold without boundary.
Suppose $u$  is a positive solution to the inequality
\[
 4 \Delta u - R u  + 2 u \ln u + \Lambda u \ge 0.
\leqno(2.10)
\]with $\Lambda \le 0$.

Given a nonnegative function $\phi \in C^\infty({\bf M})$, \ $\phi \le 1$, suppose there is
a smooth function $f$ such that $R \ge 0$ in the support of $\phi$ and that
\[
4 | \nabla f |^2 \le R - 2 \ln^+ u + \frac{|\Lambda|}{2}
\]also in the support of $\phi$.
Then
\[
\frac{\Lambda}{2} \Vert e^f  \phi u \Vert_2 \le 8 \left[
 \sup_{x \in supp \nabla \phi} e^f \sqrt{R - 2 \ln^+ u + \frac{|\Lambda|}{2}}
+ \Vert e^f \nabla \phi \Vert_\infty \right] \ \Vert u \Vert_2.
\]
\proof
\end{lemma}

The main point of the lemma is that the right hand side depends only on information  in
the support of $\nabla \phi$.

Using integration by parts,
\[
\al
&\int e^f \phi u \left( -4 \Delta + R - 2 \ln u -\Lambda - 4 |\nabla f|^2 \right) (e^f \phi u)
\\
& = 4 \int |\nabla (e^f \phi u)|^2 + \int (e^f \phi u)^2 (R - 2 \ln u -\Lambda- 4 |\nabla f|^2).
\eal
\]By assumption
\[
R-2 \ln u -\Lambda - 4 |\nabla f|^2 \ge |\Lambda|/2.
\]Hence
\[
\int e^f \phi u \left( -4 \Delta + R - 2 \ln u -\Lambda - 4 |\nabla f|^2 \right) (e^f \phi u)
\ge \frac{\Lambda}{2} \int (e^f \phi u)^2 .
\leqno(2.11)
\]By straight forward calculation
\[
\al
\text {left side of (2.11)} &= \int
(e^f \phi)^2 u \left( -4 \Delta u + R u- 2 u\ln u -\Lambda u \right) \\
&\qquad - \int e^f \phi u \left[ 8 \nabla (e^f \phi) \nabla u +
4 \Delta(e^f \phi) u \right] - 4 \int  (e^f \phi u)^2 |\nabla f|^2\\
&\le - \int e^f \phi u \left[ 8 \nabla (e^f \phi) \nabla u +
4 \Delta(e^f \phi) u \right] - 4 \int  (e^f \phi u)^2 |\nabla f|^2.
\eal
\]The last step is due to (2.10).  This together with (2.11) yield
\[
\frac{\Lambda}{2} \int (e^f \phi u)^2  \le
- \int e^f \phi u \left[ 8 \nabla (e^f \phi) \nabla u +
4 \Delta(e^f \phi) u \right] - 4 \int  (e^f \phi u)^2 |\nabla f|^2.
\]
Performing integration by parts on the term containing $\Delta$, we deduce
\[
\frac{\Lambda}{2} \int (e^f \phi u)^2  \le
- 8 \int e^f \phi u  \nabla (e^f \phi) \nabla u + \int
4 \nabla(e^f \phi) \nabla(e^f \phi u^2)  - 4 \int  (e^f \phi u)^2 |\nabla f|^2.
\]This shows
\[
\frac{\Lambda}{2} \int (e^f \phi u)^2  \le
 4 \int
| \nabla(e^f \phi)|^2 u^2 - 4 \int  (e^f \phi u)^2 |\nabla f|^2.
\]Hence
\[
\frac{\Lambda}{2} \int (e^f \phi u)^2  \le
 4 \int \left[
(e^f \phi)|^2 |\nabla f|^2 + 2 e^{2f} (\nabla f \nabla \phi )\phi + e^{2f} |\nabla \phi|^2 \right]
u^2 - 4 \int  (e^f \phi u)^2 |\nabla f|^2.
\]The first and the last term on the right hand side cancel to give
\[
\frac{\Lambda}{2} \int (e^f \phi u)^2  \le
 8 \int
  e^{2f} (\nabla f \nabla \phi )\phi u^2 + 4 \int e^{2f} |\nabla \phi|^2
u^2.
\]Note that the integrations on the right side only take place in the support
 of $\nabla \phi$. Thus it shows, by assumption on $|\nabla f|^2$,
\[
\al
&\frac{\Lambda}{2} \int (e^f \phi u)^2  \le
 4 \int_{supp \nabla \phi}
  e^{2f} |\nabla f|^2 \phi^2 u^2 + 8 \int e^{2f} |\nabla \phi|^2
u^2\\
&\le \int_{supp \nabla \phi}
  e^{2f} (R-2\ln^+ u + \frac{|\Lambda|}{2})  \phi^2 u^2 + 8 \int e^{2f} |\nabla \phi|^2
u^2.
\eal
\]So finally
\[
\frac{\Lambda}{2} \int (e^f \phi u)^2  \le \sup_{x \in supp \nabla
\phi} e^{2 f} (R-2\ln^+ u + \frac{|\Lambda|}{2})  \ \int u^2 + 8
\sup e^{2f} |\nabla \phi|^2 \int u^2.
\]\qed

\medskip

\begin{lemma}
Let $({\bf M}, g)$ be any compact manifold without boundary and ${\bf X}$ be a
domain in ${\bf M}$.
Define
\[
\lambda_X = \inf \{ \int ( 4 |\nabla v|^2 + R v^2- v^2 \ln v^2) \ | \
v \in C^\infty_0({\bf X}), \ \Vert v \Vert_2 = 1 \},
\leqno(2.12)
\]
\[
\lambda_M= \inf \{ \int ( 4 |\nabla v|^2 + R v^2- v^2 \ln v^2) \ | \
v \in C^\infty({\bf M}), \ \Vert v \Vert_2 = 1 \},
\leqno(2.13)
\]

Let $u(>0)$ be the minimizer for $\lambda_M$. For any smooth cut-off
function $\eta \in C^\infty_0({\bf X})$, \ $0 \le \eta \le 1$, it
holds
\[
\lambda_X \le \lambda_M + 4 \frac{ \int  u^2 |\nabla \eta|^2}{\int (u \eta)^2} -
\frac{\int (u \eta)^2 \ln \eta^2}{\int (u \eta)^2}.
\]
\proof
\end{lemma}

Since $\eta u/\Vert \eta u \Vert_2 \in C^\infty_0({\bf X})$  and it is $L^2$ norm is $1$,
we have, by  definition,
\[
\lambda_X \le \int \left[ 4 \frac{|\nabla (\eta u)|^2}{\Vert \eta u \Vert^2_2 } +
R \frac{(\eta u)^2}{\Vert \eta u \Vert^2_2 } - \frac{(\eta u)^2}{\Vert \eta u \Vert^2_2 }
\ln \frac{(\eta u)^2}{\Vert \eta u \Vert^2_2 } \right].
\]This implies
\[
\lambda_X \Vert \eta u \Vert^2_2  \le \int \left[ 4 |\nabla (\eta u)|^2 +R (\eta u)^2 -
 (\eta u)^2 \ln (\eta u)^2 \right] + \Vert \eta u \Vert^2_2 \ln  \Vert \eta u \Vert^2_2.
 \leqno(2.14)
 \]

 On the other hand, $u$ is a smooth positive solution (cf [Ro]) of the equation
 \[
 4 \Delta u - R u + 2 u \ln u + \lambda_M u=0.
 \]Using $\eta^2 u$ as a test function for the equation, we deduce
 \[
 \lambda_M \int (\eta u)^2 = - 4 \int (\Delta u) \eta^2 u + \int R (\eta u)^2 -
 2 \int (\eta u)^2 \ln u.
 \]By direct calculation
 \[
 -4 \int (\Delta u)  \eta^2 u = 4 \int |\nabla (\eta u)|^ 2- 4 \int u^2 |\nabla \eta|^2.
 \]Hence
 \[
  \lambda_M \int (\eta u)^2 = 4 \int |\nabla (\eta u)|^ 2- 4 \int u^2 |\nabla \eta|^2
 + \int R (\eta u)^2 -
 2 \int (\eta u)^2 \ln u.
 \leqno(2.15)
 \] Comparing (2.15) with (2.14) and noting that $\Vert \eta u \Vert_2 < 1$, we
 obtain
 \[
 \lambda_X \Vert \eta u \Vert^2_2  \le  \lambda_M \Vert \eta u \Vert^2_2 +
 4 \int |\nabla \eta|^2 u^2 - \int (\eta u)^2 \ln \eta^2.
 \]\qed

 \medskip
Now we are ready to give a
\medskip

 {\bf Proof of Theorem 1.1.}

At a given time $t$ in a Ricci flow $({\bf M}, g(t))$ and for
$\sigma>0$,
 let us define
 \[
 \lambda_{\sigma^2}(g(t)) = \inf \{ \int [ \sigma^2 ( 4 | \nabla v|^2 + R v^2) -
 v^2 \ln v^2 ] d\mu(g(t)) -
 n \ln \sigma \ | \ v \in C^\infty({\bf M}), \  \Vert v \Vert_2 =1 \}.
 \leqno(2.16)
 \]Sometimes, we refer to $\lambda_{\sigma^2}$ as the best Log Sobolev constant
 with parameter
 $\sigma$. If $t$ happen to be a surgery time, then
 $\lambda_{\sigma^2}(g(t^+))$ stands for the best Log Sobolev constant
 with parameter
 $\sigma$ for the manifold right after surgery; and
 \[
\lambda_{\sigma^2}(g(t^-)) \equiv \lim_{s \to t^-}
\lambda_{\sigma^2}(g(s)).
\]We will see in step 2 below that such limit exists.

 The main aim is to find a uniform
 lower bound for $\lambda_{\sigma^2}(g(t))$, $t \in [T_1, T_2]$, $\ \sigma
 \in (0, 1]$.

 The rest of the proof is divided into 5 steps.
 \medskip

 {\bf Step 1.  We estimate the change of $\lambda_{\sigma^2}(t)$, the best
 constant of the log Sobolev inequality,
 after one $(r, \delta)$ surgery.}
 \medskip

It will be clear that the proof below is independent of the number
of cut offs occurring in one surgery time $T$. Therefore we just
assume there is one $\epsilon$ horn and one cut off at $T$.

Let $({\bf M}, g(T^+))$ be the manifold right after the surgery and

\[
\Lambda \equiv \lambda_{\sigma^2}(g(T^+))
\]be the best constant for
this post surgery manifold, defined in (2.16).

By [Ro], there is a smooth positive function $u$ that reaches the infimum in (2.16)  and $u$
solves
\[
\sigma^2 ( 4 \Delta u - R u ) + 2 u \ln u + \Lambda u +
 n( \ln
\sigma) u =0. \leqno(2.17)
\]

After taking the scaling
\[
g_1= \sigma^{-2} g(T^+), \ R_1 = \sigma^2 R,  d_1 = \sigma^{-1} d,  \ u_1 = \sigma^{n/2} u
\]we see that $u_1$ satisfies
\[
4 \Delta_1 u _1 - R_1 u_1 + 2 u_1 \ln u_1 + \Lambda u_1 =0
\leqno(2.18)
\]and
\[
\Lambda = \inf \{ \int ( ( 4 | \nabla_{g_1} v|^2 + R_1 v^2 -v^2 \ln
v^2) d\mu(g_1)
  \ | \ v \in C^\infty({\bf M^+}), \  \Vert v \Vert_2 =1 \}.
 \leqno(2.19)
\]

Denote by $U$ the $ \sigma^{-1} h$ neighborhood of the surgery cap
${\bf C}$ under $g_1$, i.e.
\[
U = \{ x \in ({\bf M}, g_1(T^+))  \  | \  d_1(x,  {\bf C})<\sigma^{-1} h \} = \{ x \in M^+ \  | \  d(x,  {\bf C})< h \} .
\]Note that $U-{\bf C}$ is part of the $\epsilon$ tube which is unaffected by the surgery.
Therefore,  $U-{\bf C}$ is $\epsilon$ close to a portion of the standard round neck
under the scaled metric $\sigma^2 h^{-2} g_1$. Actually it is even $\delta(<\epsilon)$ close since it is part of the strong $\delta$ neck. But we do not need this fact.
Following the description at the beginning of the section, there is a longitudinal parametrization of $U-{\bf C}$, called $z$ which maps $U-{\bf C}$ to
$(-1, 0) \subset (-\epsilon^{-1}, \epsilon^{-1})$.
Let $\zeta: [-1, 0] \to [0, 1]$ be a smooth decreasing function such that $\zeta(-1)=1
$ and $\zeta(0)=0$. Then $\eta \equiv \zeta (z(x))$ maps $U-{\bf C}$ to $(0, 1)$.
We then extend $\eta$ to be a cut off function on the whole manifold by setting $\eta=1$
in ${\bf M^+}-U$ and $\eta = 0$ in ${\bf C}$.

Define
\[
\Lambda_X = \inf \{ \int ( ( 4 | \nabla_{g_1} v|^2 + R_1 v^2 -v^2
\ln v^2) d\mu(g_1)
  \ | v \in C^\infty_0({\bf M^+-{\bf C}}), \  \Vert v \Vert_2 =1 \}.
 \leqno(2.20)
\]Then it is clear that
\[
\lambda_{\sigma^2}(g(T^-)) \le \Lambda_X.
\]By Lemma 2.3,
\[
\Lambda_X \le \Lambda + 4 \frac{ \int  u^2_1 |\nabla_{g_1} \eta|^2
d\mu(g_1)}{\int (u_1 \eta)^2 d\mu(g_1)} - \frac{\int (u_1 \eta)^2
\ln \eta^2  d\mu(g_1)}{\int (u_1 \eta)^2  d\mu(g_1)}.
\]Observe that the supports of $\nabla_{g_1} \eta$ and $\eta \ln \eta$ are
in $U-{\bf C}$. Moreover
\[
|\nabla_{g_1} \eta| \le \frac{c \sigma}{h}, \qquad - \eta^2 \ln \eta^2 \le c.
\] Therefore
the above shows
\[
\lambda_{\sigma^2}(g(T^-)) \le \Lambda_X \le \Lambda + \frac{4 c
\sigma^2}{h^2} { \frac{  \int_U  u^2  d\mu(g_1)}{1-\int_U u^2
d\mu(g_1)}} + c \frac{\int_U u^2  d\mu(g_1)}{1-\int_U u^2
d\mu(g_1)}. \leqno(2.21)
\]

Recall that $\Lambda = \lambda_{\sigma^2}(g(T^+$)). So, in order to
bound it below, we need to show that $\int_U u^2_1 d\mu(g_1)$ is
small. This is where we will use Lemma 2.1 and 2.2.

Under the metric $g_1 = \sigma^{-2} g$, the capped $\epsilon$ horn
$D$ of radius $r$ under $g(T^+)$ is just a capped $\epsilon$ horn of
radius $r_1 = \sigma^{-1} r$. Using the longitudinal parametrization
$z$ of $D$ as described at the beginning the section, we can
construct a cut-off function $\phi=\phi(z(x))$ for $x \in M^+$,
which satisfies the following property.

i). $ \{ x \in {\bf M} \ | \ z(x) = 0 \}$ is the boundary of $D$.

ii). If $z \le 0$, then $\phi(z) =0$; and if $z \ge 1$, then
$\phi(z)=1$.

iii). $0 \le \phi \le 1;$, \qquad  $|\nabla_{g_1} \phi| \le
\frac{c}{r_1}$.

iv) $\phi$ is set to be zero outside of $D$ and is set to be $1$ to
the right of the set
\[
\{ x \in {\bf M}^+ \ |  \ z(x)=1\}.
\]

Notice that the support of $\nabla \phi$ is in the set where $z$ is
between $0$ and $1$. Applying Lemma 2.1 on $u_1$, which satisfies
(2.18), we know that
\[
u_1(x) \le c \max \{ \frac{1}{r_1}, 1 \}, \qquad x \in D.
\]Hence, for a negative number $\Lambda_0$ with $|\Lambda_0|$ being sufficiently large,
\[
\begin{cases}
 R_1(x) -2 \ln^+ u_1(x) + \frac{|\Lambda_0|}{2} \le c r^{-2}_1
+\frac{|\Lambda_0|}{2}, \qquad x \in supp \nabla_{g_1} \phi;\\
R_1(x) - 2\ln^+ u_1(x) + \frac{|\Lambda_0|}{2} \ge \frac{R_1(x)}{2}+
c r^{-2}_1 - c_1 \ln^+\max \{ \frac{1}{r_1}, 1 \} +
\frac{|\Lambda_0|}{2} \ge \frac{R_1(x)}{2} + \frac{|\Lambda_0|}{4} ,
\qquad x \in D.
\end{cases}
\leqno(2.22)
\]We stress that $\Lambda_0$ is independent of the size of
$r_1=\sigma/r$ which could be large or small due to the scaling
factor $\sigma$.

Recall that we aim to find a uniform lower bound for $\Lambda$. If
$\Lambda=\lambda_{\sigma^2}(g(T^+)) \ge \Lambda_0$, then we are in
good shape. So we assume throughout that $\Lambda \le \Lambda_0$.
Then, by (2.18), it holds
\[
4 \Delta_1 u _1 - R_1 u_1 + 2 u_1 \ln u_1 + \Lambda_0 u_1 \ge 0
\leqno(2.23)
\]

Motivated by Lemma 92.10 in [KL], we choose a function $f=f(x)$ as
the distance between $x$ and the set $z^{-1}(0)$ under the metric
\[
\frac{1}{4} (R_1(x) - 2\ln^+ u_1(x) + \frac{|\Lambda_0|}{2}) g_1(x),
\qquad x \in D.
\]By the first inequality in (2.22),  in the support of $\nabla_{g_1} \phi$,
\[
4 |\nabla_{g_1} f|^2 \le c r^{-2}_1 +\frac{|\Lambda_0|}{2}
\leqno(2.24)
\]and in $D$,
\[
4 |\nabla_{g_1} f|^2 \le R_1(x) - 2\ln^+ u_1(x) +
\frac{|\Lambda_0|}{2}. \leqno(2.25)
\]Note that the right hand side of (2.25) is positive by the second
inequality in (2.22).

Inequalities (2.25) and (2.23) allow us to use Lemma 2.2 (with
$\Lambda$ there replaced by $\Lambda_0$ here) to conclude
\[
\frac{\Lambda_0}{2} \Vert e^f  \phi u_1 \Vert_2 \le 8 \left[ \sup_{x
\in supp \nabla_{g_1} \phi} e^f \sqrt{R_1 - 2\ln^+ u_1 +
\frac{|\Lambda_0|}{2}} + \Vert e^f \nabla_{g_1} \phi \Vert_\infty
\right] \ \Vert u_1 \Vert_2.
\]Here the underlying  metric is $g_1$. By (2.22) (first item) this shows
\[
\frac{\Lambda_0}{2} \Vert e^f  \phi u_1 \Vert_2 \le c  \sup_{x \in
supp \nabla_{g_1} \phi} e^f \sqrt{(\frac{1}{r^2_1} + |\Lambda_0|)}
\quad
 \Vert u_1 \Vert_2.
\leqno(2.26)
\]

From (2.26), we will derive a bound for $\Vert u_1 \Vert_{L^2(U)}$
which holds for all finite $\sigma$. Here and later $\Vert u_1
\Vert_{L^2(U)}$ stands for integration under the metric $g_1$.

 First, we note from (2.26)
\[
\frac{\Lambda_0}{2} \inf_{U} e^f \Vert  u_1 \Vert_{L^2(U)} \le c
\sup_{x \in supp \nabla_{g_1} \phi} e^f \
 \sqrt{(\frac{\sigma^2}{r^2} + |\Lambda_0|)} \quad
 \Vert u_1 \Vert_2.
 \leqno(2.27)
\]

Let us remember that $U$ lies deep inside the capped $\epsilon$ horn $D$.
Going from $\partial D$ (i.e. $z^{-1}(0)$) to $U$, one must traverse a number of
disjoint $\epsilon$ necks. The ratio of scalar curvatures between the two ends of
a $\epsilon$ neck is bounded by $e^{c_2 \epsilon}$ for some fixed $c_2>0$.
The ratio of the scalar curvatures between $ \partial U$ and $\partial D$ is $c_3
r^2 h^{-2}$, which is independent of the scaling factor $\sigma$. Therefore one
must traverse a least
\[
K \equiv \frac{1}{c_2 \epsilon} \ln (c_3 r^2 h^{-2})
\]number of
$\epsilon$ necks to reach $U$. Note $K$ is independent of $\sigma$.

Let $G_i$ be one of the $\epsilon$ necks. The distance between its two ends
under the metric $g$ is comparable to $ 2 \epsilon^{-1} R^{-1/2}(x_i)$ where $x_i$
is a point in $G_i$.  So, under the metric,
\[
\frac{1}{4} (R_1(x) - 2 \ln^+ u_1(x) + \frac{|\Lambda_0|}{2}) g_1(x)
\]The distance between the two ends is bounded from below by
\[
c_4 \inf_{x \in G_i} \sqrt{\frac{1}{4} (R_1(x) - 2 \ln^+ u_1(x) +
\frac{|\Lambda_0|}{2}) } R^{-1/2}_1(x_i) \epsilon^{-1} \ge c_5
\epsilon^{-1}.
\]Here the last inequality comes from the second item in (2.22).
This means that the function $f$ increases by at least $c_5
\epsilon^{-1}$ when traversing one $\epsilon$ neck.

Next we observe that
\[
\inf_{G_2} f \ge \sup_{ supp \nabla_{g_1} \phi} f
\]since the support of $\nabla_{g_1} \phi$ is contained in the
 first $\epsilon$ neck $G_1$.
Therefore
\[
\inf_U f \ge c_5 \epsilon^{-1} (K-2) +  \inf_{G_2} f \ge c_5
\epsilon^{-1} (K-2) + \sup_{ supp \nabla_{g_1} \phi} f.
\]Substituting this to (2.27), we deduce
\[
 \Vert  u_1 \Vert_{L^2(U)} \le
2 c \Lambda^{-1}_0   e^{- c_5 \epsilon^{-1} (K-2)}   \
 \sqrt{(\frac{\sigma^2}{r^2} + |\Lambda_0|)} \quad
 \Vert u_1 \Vert_2.
\]Therefore, by the formula for $K$ in the above,
\[
 \Vert  u_1 \Vert_{L^2(U)} \le
c_6 \Lambda^{-1}_0      ( r^{-2} h^2)^{c_7 \epsilon^{-2}}    \
 \sqrt{(\frac{\sigma^2}{r^2} + |\Lambda_0|)} \quad
 \Vert u_1 \Vert_2.
\]Since $r \le$ by assumption, we know that
\[
 \Vert  u_1 \Vert_{L^2(U)} \le
c_8 C(\Lambda_0)  \sigma r^{-1}     ( r^{-2} h^2)^{c_7 \epsilon^{-2}}    \
 \Vert u \Vert_2.
\]Since $h \le \delta^2 r \le 1$, it is easy to see that we can choose $\delta$ as a suitable
power of $r$ so that
\[
\Vert  u \Vert_{L^2(U, d\mu(g))}= \Vert  u_1 \Vert_{L^2(U)} \le c_9
\sigma h^5  \Vert u \Vert_2
 \leqno(2.28)
 \]if $\epsilon$ is made sufficiently small, once and for all.

Substituting (2.28) to (2.21), we see that
\[
\lambda_{\sigma^2} (g(T^-)) \le \Lambda + c_{10} \sigma^3 h^3
\frac{1}{1-c_9 \sigma h^5}.
\]
Hence, given any $\sigma_0>0$, we have, for all $ \sigma \in (0, \sigma_0)$,
either
\[
\lambda_{\sigma^2} (g(T^+))  \ge \Lambda_0
\]or

\[
\lambda_{\sigma^2} (g(T^-)) \le \Lambda + c_{11} \sigma^3 h^3 =
\lambda_{\sigma^2} (g(T^+))  + c_{11} \sigma^3 h^3.
\]provided that $h \le ( 2 \sigma_0 c_9)^{-1/5}$.
This shows, for all $\sigma \in (0, \sigma_0]$, either
$\lambda_{\sigma^2} (g(T^+))  \ge \Lambda_0$, or
\[
\lambda_{\sigma^2} (g(T^-)) \le \lambda_{\sigma^2} (g(T^+))  +
c_{12} | vol({\bf M}(T^-))-vol({\bf M}(T^+))| \leqno(2.29)
\]
Here $vol({\bf M}(T^-))$ is the volume of the pre-surgery manifold
at $T$ and $vol({\bf M}(T^+))$ is the volume of the post-surgery
manifold at $T$.

\medskip

{\bf Step 2. We estimate the change of the best constant in the log
Sobolev inequality in a given time interval without surgery.}
\medskip

Suppose the Ricci flow is smooth from time $t_1$ to $t_2$. Let $t
\in (t_1, t_2)$ and $\sigma>0$. Recall that, for $({\bf M}, g(t))$,
Perelman's $W$ entropy with parameter $\tau$ is
\[
W(g, f, \tau)=\int_{\bf M}
                   \left(
                    \tau (  R+|\nabla f|^2)+f-n
                   \right)\tilde u\,d\mu(g(t))
\]where $\tilde u=\frac{e^{-f}}{(4 \pi \tau)^{n/2}}$. We are using
$\tilde u$ in this step to distinguish from $u$ in the last step.

We define
\[
\tau = \tau(t) = \sigma^2 + t_2 -t
\]
so that $\tau_1 = \epsilon^2+t_2-t_1$ and $\tau_2=\sigma^2$ (by
taking $t=t_1$ and $t=t_2$ respectively).

Let $\tilde u_2$ be a minimizer of the entropy $W(g(t), f, \tau_2)$
for all $\tilde u$ such that $\int \tilde{u} d\mu(g(t_2))=1$.

We solve the conjugate heat equation with the final value chosen as
$\tilde u_2$ at $t=t_2$. Let $\tilde u_1$ be the value of the
solution of the conjugate heat equation at $t=t_1$. As usual, we
define functions $f_i$ with $i=1, 2$ by the relation $\tilde u_i =
e^{-f_i}/(4 \pi \tau_i)^{n/2}$, $i=1, 2$. Then, by the monotonicity
of the $W$ entropy ([P1])
\[
 \al inf_{\int \tilde u_0 d\mu(g(t_1))=1} W(g(t_1), f_0,
\tau_1) &\le
W(g(t_1), f_1, \tau_1) \le W(g(t_2), f_2, \tau_2) \\
&= inf_{\int \tilde u d\mu(g(t_2))=1} W(g(t_2), f, \tau_2).\eal
\]
Here $f_0$ and $f$ are given by the formulas
\[
\tilde u_0 = e^{-f_0}/(4 \pi \tau_1)^{n/2}, \qquad \tilde u =
e^{-f}/(4 \pi \tau_2)^{n/2}.
\]
Using these notations we can rewrite the above as
\[
\al
& \inf_{\Vert \tilde u \Vert_1=1} \int_{\bf M}
\left(
                   \sigma^2 (R+|\nabla \ln \tilde u |^2)- \ln \tilde u
                   -\ln(4\pi \sigma^2)^{n/2}
                   \right) \tilde u\,d\mu(g(t_2))\\
&\ge \inf_{\Vert \tilde u_0 \Vert_1=1} \int_{\bf M} \left(
                   (\sigma^2+t_2-t_1) (R+|\nabla \ln \tilde u_0 |^2)- \ln
                   \tilde u_0
                   -\ln(4\pi (\sigma^2 +t_2-t_1))^{n/2}
                   \right) \tilde u_0\,d\mu(g(t_1)).
\eal
\]Denote $v=\sqrt{\tilde u}$ and $v_0 =\sqrt{\tilde u_0}$. This inequality is converted to
\[
\al
 &\inf_{\Vert v \Vert_2=1} \int_{\bf M}                \left(
                   \sigma^2 (R v^2+4 |\nabla v |^2)- v^2 \ln v^2
                   \right) \,d\mu(g(t_2)) - \ln(4\pi \sigma^2)^{n/2}\\
&\ge \inf_{\Vert v_0 \Vert_2=1} \int_{\bf M} \left(
                   4(\sigma^2+t_2-t_1) (\frac{1}{4}R v^2_0+|\nabla v_0 |^2)- v^2_0 \ln
                   v^2_0
                   \right) \,d\mu(g(t_1)) -\ln(4\pi (\sigma^2 +t_2-t_1))^{n/2}.
\eal
\]That is
\[
\lambda_{\sigma^2}(g(t_2)) \ge \lambda_{\sigma^2+t_2-t_1}(g(t_1)).
\leqno(2.30)
\]
\medskip

{\bf Step 3. We estimate the change of the best constant in the log
Sobolev inequality in the time interval $[T_1, T_2]$, with
surgeries.}
\medskip

Now, let
\[
T_1 \le t_1<t_2<...<t_k \le T_2
\]and $t_i$, $i=1, 2, ..., k$ be all the surgery times from $T_1$ to $T_2$.
Here, without loss of generality, we assume that $T_1$ and $T_2$ are not surgery
times. Otherwise we can just directly apply step 1 two more times at $T_1$ and
$T_2$.  We also fix a
\[
\sigma_0=T_2-T_1+1,
\]where $\sigma_0$ is the upper bound for the parameter $\sigma$ in step 1, (2.29).

For any $\sigma \in (0, 1]$, by (2.30), we have
\[
\lambda_{\sigma^2}(g(T_2)) \ge \lambda_{\sigma^2+T_2-t_k}(g(t^+_k)).
\] By step 1 ((2.29)),
either
\[
 \lambda_{\sigma^2+T_2-t_k}(g(t^+_k)) \ge \Lambda_0
 \]or
 \[
  \lambda_{\sigma^2+T_2-t_k}(g(t^+_k)) \ge  \lambda_{\sigma^2+T_2-t_k}(g(t^-_k))
-c_{12} |vol({\bf M}(t^-_k)-vol({\bf M}(t^+_k))|.
\]In the first case,
we have
\[
\lambda_{\sigma^2}(g(T_2)) \ge \Lambda_0.
\]So a uniform lower bound is already found.

In the second case,
\[
\lambda_{\sigma^2}(g(T_2)) \ge \lambda_{\sigma^2+T_2-t_k}(g(t^-_k))
-c_{12} |vol({\bf M}(t^-_k)-vol({\bf M}(t^+_k))|.
\]
From here we start with $\lambda_{\sigma^2+T_2-t_k}(g(t^-_k))$ and
repeat the above process. We have, from (2.30), with $\sigma^2$ in
(2.30) replaced by $\sigma^2 + T_2-t_k$,
\[
\lambda_{\sigma^2+T_2-t_k}(g(t^-_k)) \ge
\lambda_{\sigma^2+T_2-t_{k-1}}(g(t^+_{k-1})).
\]

Continue like this, until  $T_1$,  we have either
\[
\lambda_{\sigma^2}(g(T_2)) \ge \lambda_{\sigma^2+T_2-T_1}(g(T_1)) -
c_{12} \Sigma^k_{i=1}  |vol({\bf M}(t^-_i)-vol({\bf M}(t^+_i))|
\] or
\[
\lambda_{\sigma^2}(g(T_2)) \ge \Lambda_0-
c_{12} \Sigma^k_{i=1}  |vol({\bf M}(t^-_i)-vol({\bf M}(t^+_i))|.
\]Note that the above process can be carried out since all the parameters under
$\lambda$ is bounded from above by $\sigma_0$.

 It is known that
\[
\Sigma^k_{i=1}  |vol({\bf M}(t^-_i)-vol({\bf M}(t^+_i))| \le \sup_{t \in [T_1, T_2]}
vol({\bf M}(t)).
\]Hence, either
\[
  \lambda_{\sigma^2}(g(T_2)) \ge \lambda_{\sigma^2+T_2-T_1}(g(T_1)) -c_{12}
  \sup_{t \in [T_1, T_2]} vol({\bf M}(t)),
  \leqno(2.31)
  \]or
  \[
    \lambda_{\sigma^2}(g(T_2)) \ge \Lambda_0-c_{12}
  \sup_{t \in [T_1, T_2]} vol({\bf M}(t)).
  \leqno(2.32)
 \]In either case, the lower bound is independent of the number of surgeries.

If (2.31) holds, then we have to find a lower bound for
$\lambda_{\sigma^2+T_2-T_1}(g(T_1))$, which is independent of
$\sigma$. Remember that it is assumed that $({\bf M}, g(T_1))$
satisfies a Sobolev inequality with constant $A_1$.  It is well
known that this implies a log Sobolev inequality. Indeed, from
\[
\bigg( \int v^{2n/(n-2)} d\mu(g(T_1)) \bigg)^{(n-2)/n} \le A_1 \int (
4 |\nabla v |^2 +  R v^2 ) d\mu(g(T_1)) + A_1  \int v^2
d\mu(g(T_1)) ,
\]
using H\"older inequality and Jensen inequality for $\ln$, we have:

 for those $v \in W^{1, 2}({\bf M}, g(T_1))$ such that $\Vert v \Vert_2 =1$, it holds
\[
\int v^2 \ln v^2 d\mu (g(T_1)) \le \frac{n}{2} \ln
\left( A_1 \int (
4 |\nabla v |^2 +  R v^2 ) d\mu(g(T_1)) + A_1   \right).
\leqno(2.33)
\]
Recall the elementary inequality: for all $z, q >0$,
\[
\ln z \le q z - \ln q -1.
\]By (2.33), this shows
\[
\int v^2 \ln v^2 d\mu (g(T_1)) \le \frac{n}{2} q
\left( A_1 \int (
4 |\nabla v |^2 +  R v^2 ) d\mu(g(T_1)) + A_1   \right) - \frac{n}{2} \ln q -
\frac{n}{2}.
\]Take $q$ such that $ \frac{n}{2} q A_1 = \sigma^2 + T_2 -T_1$.
Since $ \sigma \le 1$, this shows, for some $B=B(A_1, T_1, T_2,
n)>0$,
\[
\al
 &\lambda_{\sigma^2+T_2-T_1}(g(T_1))\\
 &
 \equiv \inf_{\Vert v \Vert_2=1} \int [ \ (\sigma^2+T_2-T_1)(4 |\nabla v |^2 +  R v^2 ) -
  v^2 \ln v^2 ] d\mu(g(T_1)) - \frac{n}{2} \ln (\sigma^2 + T_2 -T_1)  \ge -B.
\eal
\]Therefore we can conclude from (2.31) and (2.32) that
\[
 \lambda_{\sigma^2}(g(T_2)) \ge \min\{-B, \Lambda_0\} -
c_{12}
  \sup_{t \in [T_1, T_2]} vol({\bf M}(t)) \equiv A_2
 \]for all $\sigma \in (0, 1]$.  By definition ((2.16)), this is nothing but a (restricted) log Sobolev inequality for
 $({\bf M}, g(T_2))$. i.e.
 \[
 \int v^2 \ln v^2 d\mu (g(T_2)) \le
\sigma^2 \int (
4 |\nabla v |^2 +  R v^2 ) d\mu(g(T_2))  - \frac{n}{2} \ln \sigma^2 - A_2
\leqno(2.34)
\]where $\sigma \in (0, 1]$.

\medskip

{\bf Step 4.  The log Sobolev inequality (2.34) implies certain heat
kernel estimate.}
\medskip

 Let $p(x, y, t)$ be the heat kernel of $\Delta - \frac{1}{4} R$ in
 $({\bf M}, g(T_2))$ (with the fixed metric $g(T_2))$. Then (2.34) implies,
 for $t \in (0, 1]$,
 \[
 p(x, y, t) \le \exp ( 4(T_0+1) +\frac{n}{2} \ln
                   A_2 + c + R^-_0) \frac{1}{(4 \pi t)^{n/2}} \equiv
                   \frac{\Lambda}{t^{n/2}}.
 \leqno(2.35)
 \]This follows from a generalization of Davies' argument [Da],
  as done in [Z1].
 We omit the details.

 \medskip

{\bf Step 5.  The  heat kernel estimate (2.35) implies Sobolev
inequality perturbed with scalar curvature $R$ and strong
non-collapsing.}
\medskip

 This is more or less standard.  By adapting the standard method in heat
  kernel estimate in [Da],  as demonstrated in the paper [Z1],
it is known that (2.35) implies the desired Sobolev imbedding for
$g(T_2)$, i.e. for all $v \in W^{1, 2}({\bf M}, g(T_2))$.
\[
\bigg( \int v^{2n/(n-2)} d\mu(g(T_2)) \bigg)^{(n-2)/n} \le A_2 \int (
4 |\nabla v |^2 + R v^2 ) d\mu(g(T_2)) + A_2  \int v^2
d\mu(g(T_2)).
\]This is the desired Sobolev inequality.

The strong noncollapsing result follows from the work of Carron
[Ca], as given in [Z2]. Please see Lemma A.2 in the appendix. \qed

  \section{Appendix}

We collect some basic facts concerning Ricci flow with surgery. For details, please
see Perelman's papers [P1, 2] and [CZ], [KL] and [MT].
\medskip

{\bf Definition.} {\bf $(r, \delta)$} surgery.

A surgery occurs at a
$\delta$ neck, called ${\bf N}$, of radius $h$ such that $({\bf N},
h^{-2} g)$ is $\delta$ close in the $C^{[\delta^{-1}]}$ topology to
the standard round neck $S^2 \times (-\delta^{-1}, \delta^{-1})$ of
scalar curvature $1$. Let $\Pi$ be the diffeomorphism from the
standard round neck to ${\bf N}$ in the definition of $\delta$
closeness. Denote by $z$ for a number in $(-\delta^{-1},
\delta^{-1})$.  For $\theta \in S^2$, $(\theta, z)$ is a
parametrization of ${\bf N}$ via the diffeomorphism $\Pi$. We can
identify the metric on ${\bf N}$ with its pull back on the round
neck by $\Pi$ in this manner.

Following  the notations on p424 of [CZ] (based on [Ha]), the metric $\tilde
g=\tilde g(T_2)$ right after the surgery is given by
\[
\tilde g =
\begin{cases}
\bar g, \qquad z \le 0,\\
e^{-2 f} \bar g, \qquad z \in [0, 2],\\
\phi e^{-2 f} \bar g + (1-\phi)e^{-2 f} h^2 g_0 , \qquad z \in [2, 3],\\
e^{-2 f} h^2 g_0 , \qquad z \in [3, 4].
\end{cases}
\]Here $\bar g$ is the nonsingular part of the  $\lim_{t \to
T^-_2} g(t)$; $g_0$ is the standard metric on the round neck; and
$f=f(z)$ is a smooth function given by (c.f. p424 [CZ])
\[
\al
 &f(z)=0, \, z \le 0; \quad f(z) = ce^{-P/z}, \, z \in (0, 3]; \,
f''(z)>0, \, z \in [3, 3.9]; \\
&\quad f(z)=-\frac{1}{2} \ln (16-z^2), \, z \in [3.9, 4]. \eal
\]Here a small $c>0$ and a large $P>0$ are suitably chosen to
ensure that the Hamilton-Ivey pinching condition remains valid.
$\phi$ is a smooth bump function with $\phi=1$ for $z \le 2$ and
$\phi=0$ for $z \ge 3$.

 \medskip

 The next result (Lemma A.2 in [Z2]) relates the Sobolev imbedding to local noncollapsing of volume of geodesic balls. We
follow the idea in [Ca].
\medskip

{\it {\bf Lemma A.2.} Let $({\bf M}, g)$ be a Riemannian manifold.
Given $x_0 \in {\bf M}$ and $r \in (0, 1]$. Let $B(x_0, r)$ be a
proper geodesic ball, i.e. ${\bf M} - B(x_0, r)$ is non empty.
Suppose the scalar curvature $R$ satisfies $|R(x)| \le 1/r^2$ in
$B(x_0, r)$ and the following Sobolev imbedding holds:
 for all $v \in W^{1, 2}_0(B(x_0, r))$, and a constant $A \ge 1$,
 \[
\bigg( \int v^{2n/(n-2)} d\mu(g) \bigg)^{(n-2)/n} \le A \int (
|\nabla v |^2 + \frac{1}{4} R v^2 ) d\mu(g) + A  \int v^2 d\mu(g).
\]Then $|B(x_0, r)| \ge 2^{-(n+5) n/2} A^{-n/2} r^n.$
}
\medskip

\noindent {\bf Proof.} Since $R \le 1/r^2$, $r \le 1$ and $A \ge 1$
by assumption, the Sobolev imbedding can be simplified to
\[
\bigg( \int v^{2n/(n-2)} d\mu(g) \bigg)^{(n-2)/n} \le A \int |\nabla
v |^2  d\mu(g) + \frac{2 A}{r^2}  \int v^2 d\mu(g).
\]

Under the scaled metric $g_1 = g/r^2$, we have, for all $v \in W^{1,
2}_0(B(x_0, 1, g_1))$,
\[
\bigg( \int v^{2n/(n-2)} d\mu(g_1) \bigg)^{(n-2)/n} \le A \int
|\nabla v |^2  d\mu(g_1) + 2 A  \int v^2 d\mu(g_1).
\]Now, by [Ca] (see p33, line 4 of [H2]), it holds
\[
|B(x_0, 1, g_1)|_{g_1} \ge \min \{\frac{1}{2 \sqrt{2A}},
\frac{1}{2^{(n+4)/2} \sqrt{2A}} \}^n.
\]Therefore
\[
|B(x_0, r,, g)|_g \ge 2^{-(n+5) n/2} A^{-n/2} r^n.
\]\qed

\noindent e-mail:  qizhang@math.ucr.edu
\end{document}